\documentclass[12pt]{article}
\usepackage[margin=2.5cm]{geometry}
\usepackage{amsmath,amsthm,amsfonts,amssymb}
\usepackage[mathscr]{euscript}
\usepackage{amsmath}
\usepackage[table]{xcolor}

\usepackage{graphicx,latexsym}
\usepackage{lscape}
\usepackage{fixmath}
\usepackage{multicol}
\usepackage{graphicx}
\usepackage{caption}
\usepackage{float}
\usepackage{cite}
\usepackage{subfig}
\usepackage{setspace}
\usepackage{xcolor}
\usepackage{multirow}
\usepackage{hhline} 
\usepackage[utf8]{inputenc}
\usepackage{gensymb}
\usepackage{setspace}
\usepackage{abstract}
\usepackage{caption}
\newtheorem{theorem}{Theorem}[section]
\newtheorem{lemma}[theorem]{Lemma}
\newtheorem{proposition}[theorem]{Proposition}

\newcommand{\gp}{{\rm gp}}
\newcommand{\ic}{{\rm ic}}
\newcommand{\ip}{{\rm ip}}

\begin{document}
\title{General Position Problem of Butterfly Networks}

\author{
	R. Prabha$^{a}$, 
	S. Renukaa Devi$^{b,c}$,
	Paul Manuel$^{d}$
}

\date{}

\maketitle
\vspace{-0.8 cm}
\begin{center}
	$^a$ Department of Mathematics, Ethiraj College for Women, Chennai, Tamilnadu, India \\
prabha75@gmail.com\\
	\medskip

    $^b$ Research Scholar, University of Madras, Chennai, Tamilnadu, India\\ 
    \medskip
    
	$^c$ Department of Mathematics, Avichi College of Arts and Science, Chennai, Tamilnadu, India \\
	renukaadevim@gmail.com\\
	\medskip
	$^d$ Department of Information Science, College of Life Sciences, Kuwait University, Kuwait \\
	pauldmanuel@gmail.com, p.manuel@ku.edu.kw\\
	\medskip
\end{center}

\begin{abstract} 
	A \textit{general position set S}  is a set $S \subseteq V(G)$ such that no three vertices of $S$ lie on a shortest path in \textit{G}. Such a set of maximum size in \textit{G} is called a \textit{\gp-set} of \textit{G} and its cardinality is called the \textit{\gp-number} of \textit{G} denoted by \textit{\gp(G)}. The authors who introduced the general position problem stated that the general position problem for butterfly networks was open \cite{MaKl18b}. A well-known technique to solve  the general position problem for a given network is to use its isometric path cover number as an upper bound.  The general position problem for butterfly networks remained open because this technique is not applicable for butterfly networks. In this paper, we adopt a new technique which uses the isometric cycle cover number as its upper bound. This technique is interesting and useful because it opens new avenues to solve the general position problem for networks which do not have solutions yet.
\end{abstract}


\medskip
\noindent{\bf AMS Subj.\ Class.~(2020)}: 05C12, 05C82.

\section{Introduction}

For the definitions and terminologies refer \cite{BoMu08}. The $distance$ $d(x, y)$ between any two vertices $x, y$ in a graph $G$ is the length of a shortest $x, y$-path in $G$; any such path is called a $geodesic$. A subgraph, $H$ in $G$ such that $d_G(x, y) = d_H(x, y)$ for all $x, y\in V(H)$ is called an $isometric$ $subgraph$. The minimum number of isometric cycles that cover the vertices of $G$ is referred as an $isometric$ $cycle$ $cover$ $number$, $ic(G)$. 

A set $S \subseteq V(G)$ is defined as a $general$ $position$ $set$ if no three vertices lie on a shortest path in $G$. A $\gp$-$set$ of $G$ is a general position set of maximum cardinality in $G$ which is denoted by $\gp(G)$. If $G$ is a graph, the general position problem is to find its largest size general position set of $G$. In \cite{MaKl18} Paul Manuel et.al. introduced  the general position problem motivated by the no-three-in-line and general position subset selection problems \cite{Du17, FrKaNicNie17, PaWo13} and further proved it to be NP-complete. The above problem has been researched in many articles \cite{AnUlChKlTh19, GhKlMaMoRaRu19, KlPaRuYe19, KlYe19, ThCh20}. In \cite{MaKl18b} Paul Manuel et.al. has determined $\gp(G)$ of infinite grids and infinite diagonal grids using the strategy of Monotone Geodesic Labelling. Further the general position problem for Bene\v{s} networks by using isometric path covers has been solved. In the same paper, the authors have commented that the strategy of isometric path covers could not be applied for butterfly networks and claim that it remains a challenge to prove that gp-number of $r$-dim butterfly is $2^r$. In this paper, we use a novel technique by considering the isometric cycle cover of butterfly networks and hence solve the general position problem of butterfly networks.

\section{Butterfly Network}

\begin{figure}[H]
	\centering
	\includegraphics[width=1.1\linewidth]{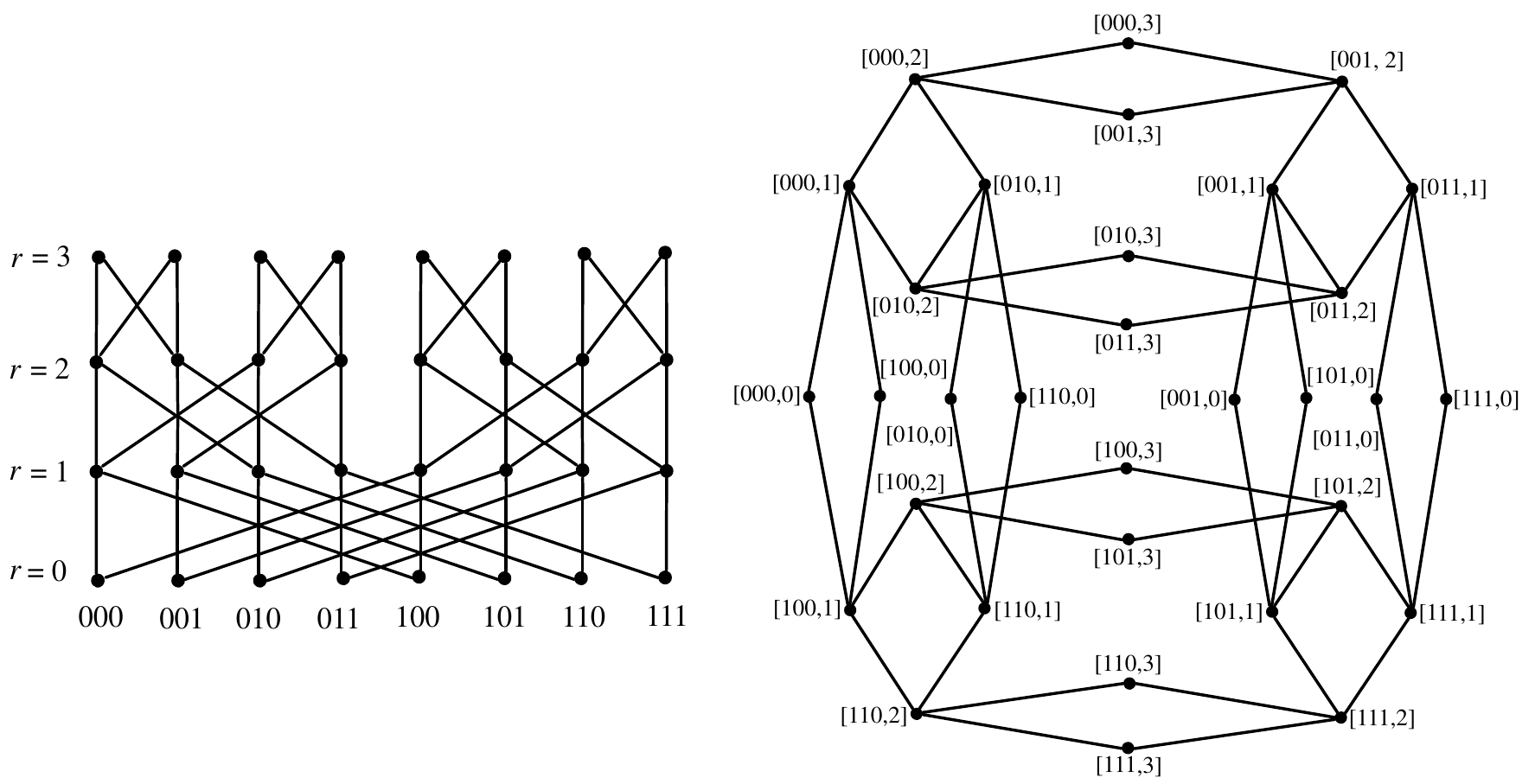}
	\caption{Normal representation and diamond representation of $BF(3)$}
	\label{bf3}
\end{figure}

\begin{figure}[H]
	\centering
	\includegraphics[width=0.6\linewidth]{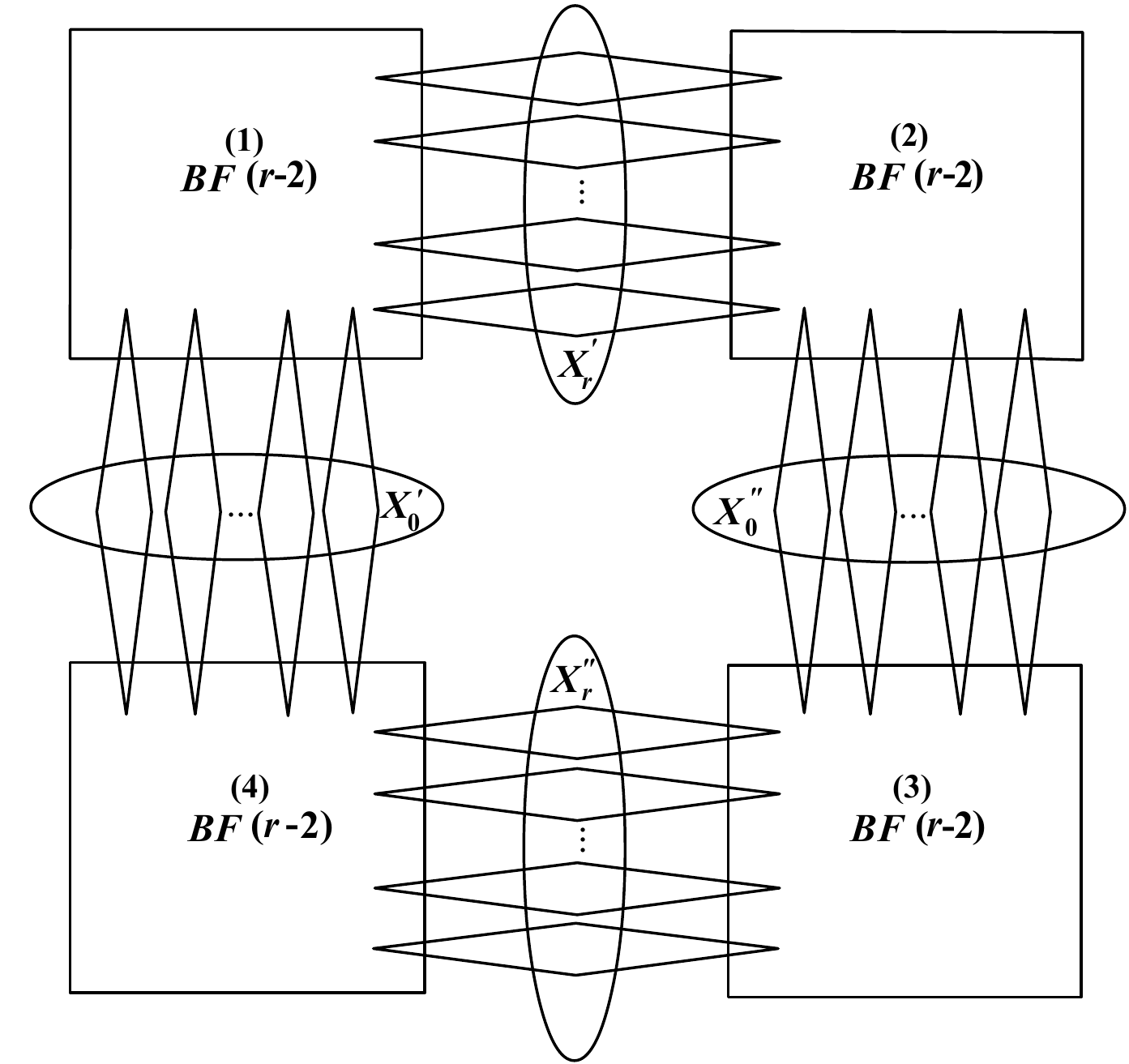}
	\caption{$BF(r)$}
	\label{bfr-2}
\end{figure}

An $r$-$dim$ $butterfly$ $network$ $BF(r)$ has vertices $[w, q]$, where $q\in \{0, 1\}^{r}$ and $w \in \{0, 1, \dots , r\}$. The vertices $[w, q]$ and $[w', q']$ are adjacent if $|w - w'|=1$, and either $q = q'$ or $q$ and $q'$ differ precisely in the $w^{th}$ bit. $BF(r)$ has $r+1$ levels with $2^r$ vertices at each level and $r2^{r+1}$ edges. The vertices at level 0 and $r$ are of 2-degree vertices and the rest are of 4-degree vertices \cite{MaAbraRa08, RaRaVe09}. $BF(r)$ has two standard graphical representations namely normal and diamond representation, for further details of which one may refer \cite{MaAbraRa08}. We shall use the diamond representation of $BF(r)$ throughout this paper. For our convenience, we denote the set of all 2-degree and 4-degree vertices by $X$ and $Y$ respectively. Further we denote any vertex $[w, q]$ of $BF(r)$ by $[a_1 a_2 \cdots a_r, r], a_i \in \{0, 1\}$ (Refer Fig. \ref{bf3}). Observe that $BF(r)$ has four disjoint copies of $BF(r-2)$ which we denote as $BF^{(1)}(r-2), BF^{(2)}(r-2), BF^{(3)}(r-2), BF^{(4)}(r-2)$. The vertices of $X$ at level 0 and level $r$ are denoted as $X_0$ and $X_r$ respectively. Let $X_0^{'}$, $X_0^{''}$, $X_r^{'}$ and $X_r^{''}$ denote the vertices of $X$ joining $BF^{(1)}(r-2)$ and $BF^{(4)}(r-2)$, $BF^{(2)}(r-2)$ and $BF^{(3)}(r-2)$, $BF^{(1)}(r-2)$ and $BF^{(2)}(r-2)$, $BF^{(3)}(r-2)$ and $BF^{(4)}(r-2)$ respectively (Refer Fig. \ref{bfr-2}). 

\par
\begin{figure}[H]
 	\centering
 	\includegraphics[width=0.77\linewidth]{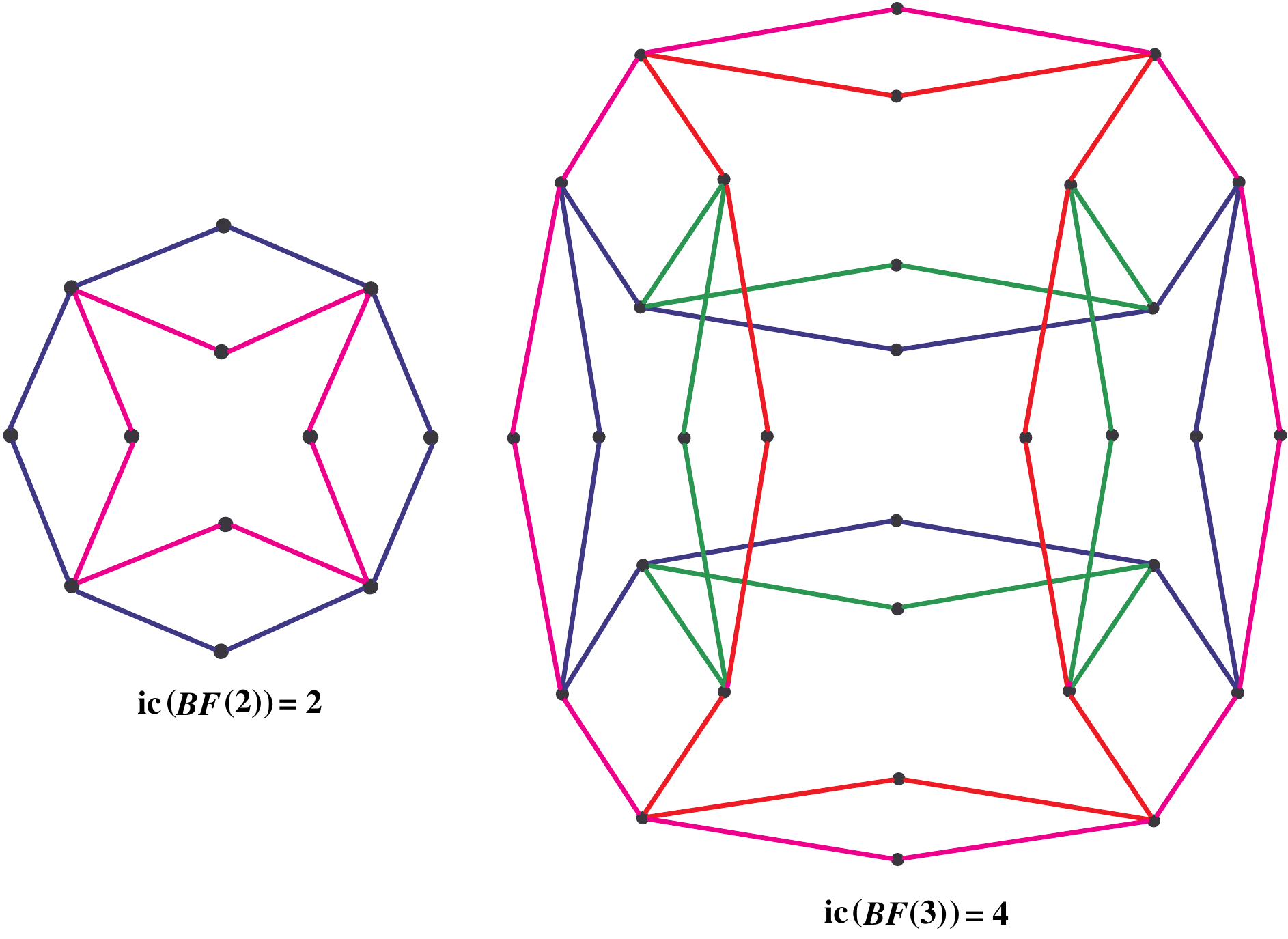}
 	\caption{Isometric cycle cover of $BF(2)$ and $BF(3)$.}
\label{icbf3}
 \end{figure}

 We now proceed as follows: first we consider the edge-disjoint isometric cycle cover of butterfly networks, constructed by the authors in \cite{PaSaPrAn22} (Refer Fig. \ref{icbf3}). Using the isometric cycle cover number as an upper bound, we compute the \gp-number of $BF(r)$. 

\par We begin by stating the results which we need from \cite{PaSaPrAn22} and \cite{MaKl18}.

\begin{lemma}\normalfont{\cite{PaSaPrAn22}}
\label{lem1-PaSaPrAn22}
\em{If $r\geq 3$, then $E(BF(r))$ can be partitioned by a set $S(r)$ of edge-disjoint isometric cycles of length $4r$, where $|S(r)| = 2^{r-1}$ and each isometric cycle of $S(r)$ has two vertices at level 0} (Refer Fig. \ref{icbf3}).
\end{lemma}

\begin{proposition} \normalfont{\cite{MaKl18}}
 	\label{prop1-MaKl18}
 	\em{Given a graph $G$, $\gp(G)\leq 2 \ip(G)$, and $\gp(G)\leq 3 \ic(G)$.}
 	 \end{proposition}

\par Next we proceed to compute an upper bound for the maximum number of 2-degree vertices in any general position set of $BF(r)$.

\begin{lemma} 
	\label{lem2-BF}
	If $S$ is a general position set of $BF(r)$ and $X$ is the set of all 2-degree vertices of $BF(r)$, then $\lvert S\cap X\rvert \leq 2^r$.
\end{lemma}

\proof
Suppose $\lvert S\cap X\rvert >2^r$. Then $S\cap X_0 \neq \emptyset$ and $S\cap X_r \neq \emptyset$. Also either $\lvert S\cap X_0\rvert > 2^{r-1}$ or $\lvert S\cap X_r\rvert > 2^{r-1}$. Without loss of generality we may suppose that $\lvert S\cap X_0\rvert > 2^{r-1}$. This implies $S\cap X_0^{'} \neq \emptyset$ and $S\cap X_0^{''} \neq \emptyset$. Choose two vertices $x, y$ such that $x\in S\cap X_0^{'}$, $y\in S\cap X_0^{''}$. Then for any $z\in S\cap X_r$, observe that $\{x, y, z\}$ is not in general position which contradicts that $S$ is a general position set of $BF(r)$. Hence the proof.

\begin{theorem}
	\label{thm1-BF}
	$\gp(BF(r)) = 2^r + 2^{r-2}$, $r\geq 2$.
\end{theorem} 

\begin{figure}[H]
	\centering
	\includegraphics[width=0.8\linewidth]{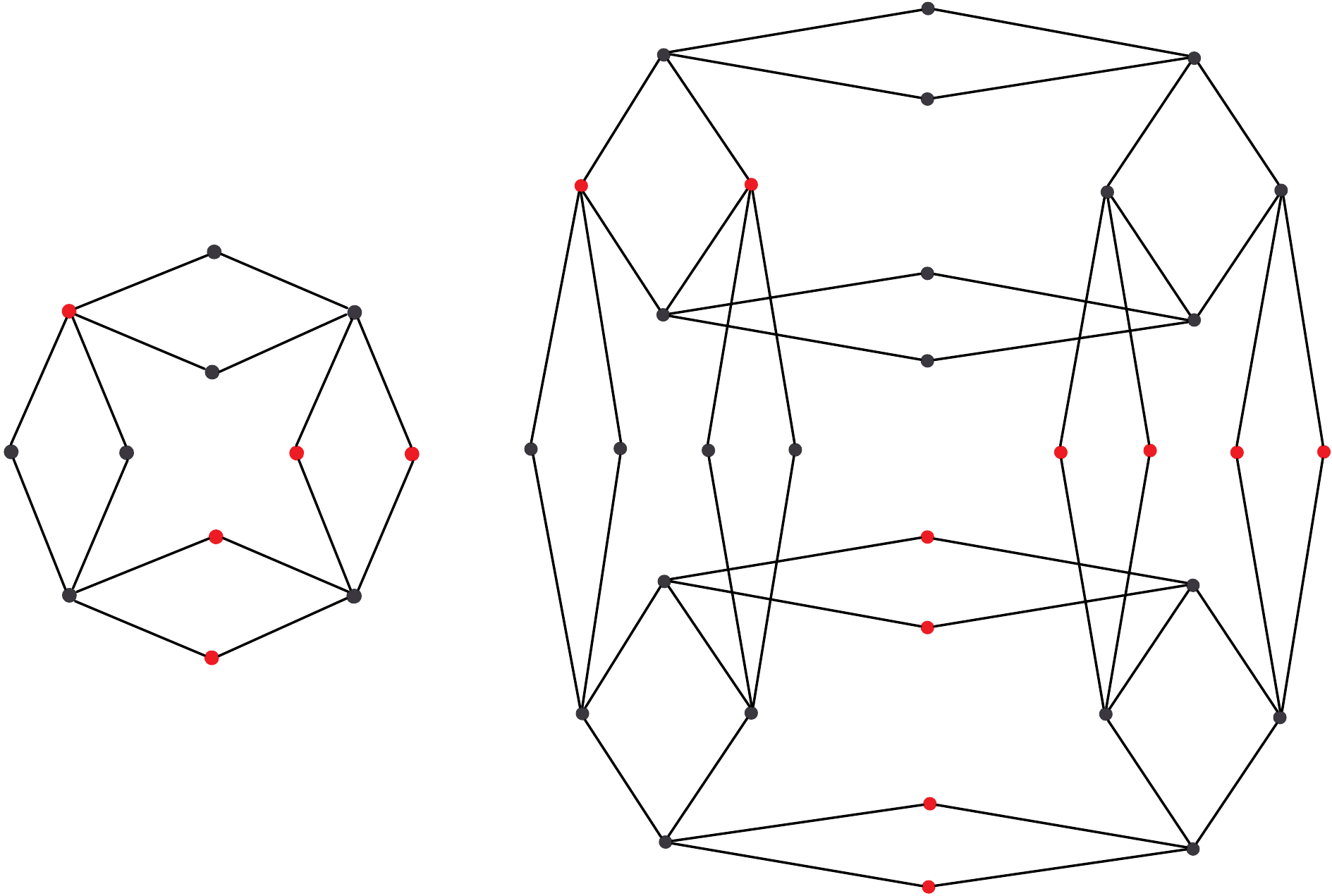}
	\caption{A gp-set of $BF(2)$ and $BF(3)$}.
	\label{bf3gp}
\end{figure}

\proof
$\gp(BF(r)) \geq 2^r + 2^{r-2}$ follows from the fact that $S = \{[a_1a_2\dots a_{r-1}1, 0]: a_1a_2\dots a_{r-1}$ is any binary sequence$\}\cup \{[1a_2a_3\dots a_r, r]: a_2a_3\dots a_r$ is any binary sequence$\}\cup \{[0a_2a_3\dots a_{r-1}0, 1]: a_2a_3\dots a_{r-1}$ is any binary sequence$\}$ is a general position set of $BF(r)$ (Refer Fig. \ref{bf3gp}).

\noindent \textbf{Claim:} $\gp(BF(r)) \leq 2^r + 2^{r-2}$.

\par Let $C_1, C_2,\dots , C_{2^{r-1}}$ denote the $2^{r-1}$ isometric cycles of $BF(r)$ as constructed in  lemma \ref{lem1-PaSaPrAn22}. Observe that every 2-degree vertex is covered by exactly one isometric cycle and every 4-degree vertex is covered by exactly two isometric cycles in $BF(r)$. In other words, exactly one isometric cycle passes through every $x\in S\cap X$ and exactly two isometric cycles pass through every $y\in S\cap Y$ (Refer Fig. \ref{icbf3}). 

\par Now let us compute $S$ under the following cases:
\\ \textbf{Case (i)} $S\cap X =\emptyset$. Then by proposition \ref{prop1-MaKl18} $\lvert S\cap Y\rvert \leq 3.2^{r-2} < 2^r$.
\\ \textbf{Case (ii)} $S\cap Y = \emptyset$. Then by lemma \ref{lem2-BF} $\lvert S\cap X\rvert \leq 2^r$.
\\ \textbf{Case (iii)} $S\cap X \neq \emptyset$ and $S\cap Y \neq \emptyset$. Now 
\\ \textbf{Subcase (i)} $\lvert S\cap (X\cap V(C_i))\rvert = 1, 1\leq i\leq 2^{r-1}$.
\par Since $\gp(C_n) = 3$, $n\geq 5$, $\lvert S\cap (Y\cap V(C_i))\rvert \leq 2, 1\leq i\leq 2^{r-2}$. Then by  proposition \ref{prop1-MaKl18} $\lvert S\rvert \leq 2^{r-1}+2.2^{r-2} \leq 2^r$.
\\ \textbf{Subcase (ii)} $\lvert S\cap (X\cap V(C_i))\rvert = 2, 1\leq i\leq 2^{r-1}$. 
\par By a similar argument $\lvert S\cap (Y\cap V(C_i))\rvert \leq 1, 1\leq i\leq 2^{r-2}$ and hence $\lvert S\rvert \leq 2^r + 2^{r-2}$. This completes the proof.

\end{document}